\documentclass[english,russian,12pt,twoside]{article}

\usepackage{amssymb,amsmath}

\usepackage[russian]{babel}

\usepackage[cp866]{inputenc}

\begin{document}

\begin{center} {\bf   Nonstationary boundary value problems  \\
for wave equation and their generalized solutions} \vspace{3mm}

{\bf Alexeyeva L.A.}
\end{center}

\centerline{\small\sc\textit{Institute of Mathematics, Kazakhstan}}

\centerline {\textit{Pushkin str.125, Almaty, 050010,
\,  alexeeva@math.kz} }
\vspace{3mm}

Investigation  of  many problems of acoustics, hydromechanics,
elastodynamics and  others  is connec\-ted with boundary value
problem (BVP) for wave equation - the multidimensional analogue of
the Dalamber's equation, which describes processes of spreading
the waves in uniform isotropic ambience. So elaboration of the
efficient methods of their solving for areas with free geometry
and varied type of the border conditions is very currently .

The most efficient method of the study of such problems is the
boundary integral equations method (BIEM). Its main advance is
reduction  of dimension of  solved equations and increasing
calculation stability. That is very greatly, particularly for
unbounded areas. The existing methods and programs of splines
approximation for free contours and surfaces remove the problem of
restriction of the form of considered areas by using BIEM. At
present time this method is broadly used for solving  elliptical
and parabolic problems, what  is connected with success in
development of BIE theories for elliptical and parabolic equations
and systems.

Solving nonstationary dynamic problems on base of the BIE method
requires entering the notion of the \textit{generalized solution}.
That is connected with particularity of the fundamental solution
of hyperbolic equations, which belongs to the class of generalized
functions. Besides classical notion to differentiability of the
decisions for hyperbolic equations sharply narrows the class of
useful for applications problems. In particular, typical physical
processes, being accompanied shock waves, are not described by
differentiated solutions of such equations.

Here  BIE method is elaborated for determination of nonstationary
solutions of wave equations in bounden region in $R^N$ with
boundary Dirihlet or Neumann conditions. The generalized solution
of BVP are constructed on the base of dynamic analogues of Green
and Gauss formulas for solutions of the wave equation in space of
generalized function. Their regular integral
re\-pre\-sen\-ta\-tions and singular BIE are built for N=1,2,3,
also at presence of shock waves. \vspace{3mm}

\textbf{1. The generalized solutions of the wave equation.Shock
waves.}The  multidimensional analogue of the Dalamber's equation
is considered:
\begin{equation}\label{(1.1)}
\Box_c u \equiv \Delta u - \frac{1}{c^ 2} \frac{{\partial ^2
u}}{{\partial t^2 }} = G(x,t),\quad x \in R^N ,\,\;t \in R^1 .
\end{equation}
$G$ is regular function, here and hereinafter $\Delta $ is Laplace
operator,  $u,_i=\frac{\partial u}{\partial
x_i},\,\,\dot{u}=u,_t=\frac{\partial u}{\partial t}$.

It is well known [1,2] that Eq. (\ref{(1.1)}) is strictly
hyperbolic, class of its solutions contains  functions with
breakup of derivatives. The surface of the breakup $F$ in $R^{N +
1} = \{(x,\tau \equiv ct )\}$ is the characteristic surface of
Eq.(\ref{(1.1)}), which satisfies to  the relations: $ \nu _\tau
^2  - \left\| \nu \right\|_N ^2  = 0,\quad \nu _\tau =\nu _{N + 1}
< 0$, where $\nu (x,\tau ) = \left( {\nu _1 ,...,\nu _N, \nu _\tau
} \right)$ is a  normal vector  to $F$ , $\left\| \nu \right\|_N =
\sqrt {\nu _j \nu _j }$.  In $R^N$ this is a wave front $F_t$,
moving  with the constant speed $c$ :
\begin{equation}\label{(1.3)}
c =  - \nu _t /\left\| \nu  \right\|_N,\quad \nu _t =c\nu _\tau
\end{equation}
On reiterative indexes $i,j$ in product everywhere their is
summation from 1 to $N$.

 There are the Adamar's conditions on the  jumps  at
$F_t$:
\begin{equation}\label{(1.4)}
\left[ {u\left( {x,t} \right)} \right]_{F_t }  = 0,\quad \left[
{\dot u + cn_i \,u,_i } \right]_{F_t }  = 0.
\end{equation}
Here $\left[ {f\left( {x,t} \right)} \right]_{F_t }  = f^ + \left(
{x,t} \right) - f^ -  \left( {x,t} \right) = \mathop {\lim
}\limits_{\varepsilon  \to  + 0} \left( {f\left( {x + \varepsilon
n,t} \right) - f\left( {x - \varepsilon n,t} \right)}
\right),\quad x \in   F_ t, $ $ n(x,t)$ is a unit normal vector to
$ F_ t $, directed aside its spreading:
\begin{equation}\label{(1.5)}
n_i  = \frac{{\nu _i }}{{\left\| \nu  \right\|_N }} =
\frac{{grad\,F }}{{\left\| {grad\,F } \right\|}},\quad i =
1,...,N;\;
\end{equation}
Last equality may be written only  if the equation of wave front
can be  present as $ F (x,t) = 0 $ , at condition of existence
$grad\,F $.

The class of the similar solutions of the hyperbolic equations is
named \textit{shock waves}.

From second condition (\ref{(1.4)}) follows
\begin{equation}\label{(1.4*)}
\dot u^ -   + cn_i \,u,_i^ -   = \dot u^ +   + cn_i \,u,_i^
+,\quad x\in F_t.
\end{equation}
If before the front  $u \equiv 0 $ , this equality gives the
useful correlation on $F_t$: $ (grad\,u,n) = - c^{ - 1} \dot u.$

We notice, that tangent derivatives to characteristic surface, on
the strength of continuity $u$, are also continues. Then
\begin{equation}\label{(1.6)}
\left[ {u,_\tau  \gamma _\tau   + u_j \gamma _j } \right]_F  =
0\quad \textrm{ for }\,\,\forall \gamma \in R^{N+1}:\,\,(\nu
,\gamma ) = 0.
\end{equation}
In particular if $\gamma = \gamma ^j = ( - \nu _j ,\,\,\nu
_\tau \delta _1^j ,\,\,\nu _\tau  \delta _2^j ,...,\,\,\nu _\tau
\delta _N^j )$, where $\delta_i^j$ is  Kronekker's symbol, then we have the
  condition of the type:
\begin{equation}\label{(1.6*)}
\left[ { - u,_\tau  \nu _j  + u_j \nu _\tau  } \right]_F  = 0
\Rightarrow \;\left[ {\dot un_j  + cu_j } \right]_{F_t }  =
0,\quad j=1,...,N.
\end{equation}

Hereinafter we shall consider the functions $u(x,t)$, which are
continues together with derivati\-ves of first and second order
almost everywhere with the exclusion of   finite or counting
number of wave fronts, on which  the conditions on jumps are
satisfied (\ref{(1.3)}). We shall name such solutions
\textit{classical}. Let's show that they are generalized solutions
of the Eq. (\ref{(1.1)}).

For this let's consider the Eq.(\ref{(1.1)}) on the space of
generalized functions  $ D'(R^{N + 1} ) = \\=\left\{ {\hat
f(x,\tau )} \right\}$, which are determined on the space of
infinitely differentiable finite functions $D(R^{N + 1} ) =
\left\{ {\varphi (x,\tau )} \right\}$  [2]. The value $ \hat f$ on
$ \varphi $, as it's accepted , is denoted as $(\hat f,\varphi )$.
For regular function $\hat f$, corresponding to local integrable
$f$, $ (\hat f,\varphi ) = \int\limits_{R^{N + 1} } {f(x,\tau )}
\varphi (x,\tau )dV(x) d\tau $.  Here everywhere $ dV(x) =
dx_1...dx_N $.

D e f i n i t i o n . Function $ \hat f \in D'(R^N )$ is
identified as \textit{generalized solution} of Eq. (\ref{(1.1)})
if for any $ \varphi \in D(R^{N + 1} )$  $ (\Box_c \hat f,\varphi
) \equiv (\hat f,\Box_c \varphi ) = (G,\varphi ).$

L e m m a 1.1. \textit{If $u(x,t)$ is  the classical solution of
Eq. (\ref{(1.1)}), then $ \hat u(x,t)$ is its generalized solution
.}

P r o o f. If $u(x,t)$ has a finite breakup on $F$ then [2] $\hat
u,_j = u,_j  + {\rm [}u{\rm ]}_F \nu _j \delta _F {\rm (}x,\tau
{\rm )}, $ where first summand on the right is a classical
derivative over $x_j ,\,\, \left\| {\bf \nu } \right\| =
1,\;\delta _F $ is \textit{simple layer} on $F$:
\[
\left( {{\rm [}u{\rm ]}_F \nu _j \delta _F {\rm (}x,\tau {\rm
)}{\rm ,}\varphi {\rm (x}{\rm ,}\tau {\rm )}} \right) =
\int\limits_F {{\rm [}u{\rm (}x{\rm ,}\tau {\rm )]}_F {\nu _j
(x,\tau)}} \varphi { (x ,\tau  )}dF(x,\tau), \quad\forall \varphi
\in D(R^{N + 1} ).
\]
Here integral on $ F$ surfaced.   On the strength of continuity
$u$ outside of wavefront\\ $ \left[ u \right]_F = \mathop {\lim
}\limits_{\varepsilon  \to  + 0} \left( {u(x + \varepsilon \,n,t)
- u(x - \varepsilon \,n,t)} \right) = u^ + (x,t) - u^ - (x,t) =
\left[ u \right]_{F_t }.$

So, with regard for (\ref{(1.5)}), we
 get: \\$\hat u,_j = u,_j + {\rm [}u{\rm ]}_{F_t } \nu _j
\delta _F {\rm (}x, \tau {\rm )},\,\,\hat u,_{jj}= u,_{jj}  +
\left[ {u,_j } \right]_{F_t } \left\| \nu \right\|_N n_j \delta _F
+\partial _j \left\{ \left\| \nu \right\|_N {\rm [}u{\rm ]}_{F_t }
n_j \delta _F \right\}.$ \\On virtue of Eq.(\ref{(1.3)})\\ $ \hat
u,_\tau = u,_\tau + {\rm [}u{\rm ]}_{F_t } \nu _\tau \delta _F =
c^{{\rm - 1}} u,_t - \left\| \nu  \right\|_N {\rm [}u{\rm ]}_{F_t
} \delta _F $, $\hat u,_{tt}  = c^{{\rm  - 2}} u,_{tt}  - c^{{\rm
- 1}} \left[ {u,_t } \right]_{F_t } \left\| \nu \right\|_N \delta
_F  - \left\{ {\left\| \nu \right\|_N {\rm [}u{\rm ]}_{F_t }
\delta _F  }\right\} _\tau .$ \\ With regard for these equalities
and Adamar's conditions (\ref{(1.1)})
\[
\displaylines{ \Box  _c \hat u =\Box _c u + \left\{ {c^{{\rm  -
1}} \left[ {u,_t } \right]_{F_t }  + \left[ {n_j u,_j }
\right]_{F_t } } \right\}\left\| \nu  \right\|_N \delta _F {\rm
(}x,\tau {\rm )} +  \cr + c^{ - 1} \partial _t \left\{ {\left\|
\nu \right\|_N {\rm [}u{\rm ]}_{F_t } \delta _F {\rm (}x,\tau {\rm
)}} \right\} +
\partial _j \left\{ {\left\| \nu  \right\|_N {\rm [}u{\rm ]}_{F_t
} n_j \delta _F {\rm (}x,\tau {\rm )}} \right\} = \hat G(x,t), }
\]
since all densities of simple and double layers on $F_t $ are
equal to zero. Really, the second summand is a zero on the
strength of the second condition (\ref{(1.4)}) on front. But two
others are null on the strength of the first  one, since their
action on $D(R^{N+1})$ is defined as
 $$
\displaylines{ \left( {c^{ - 1} \partial _t \left\{ {\left\| \nu
\right\|_N {\rm [}u{\rm ]}_{F_t } \delta _F {\rm (}x,\tau {\rm )}}
\right\} +
\partial _j \left\{ {\left\| \nu  \right\|_N {\rm [}u{\rm ]}_{F_t
} n_j \delta _F {\rm (}x,\tau {\rm )}} \right\},\varphi (x,t)}
\right) =  \cr
 = -\int\limits_F {\nu_\tau \,[u]}_F\left(\frac{{\partial \varphi
}}{{\partial n}} - \varphi ,_\tau   \right)  dF(x,\tau)=0. }$$
 The lemma has been proved.

\textit{R e m a r k \,\,  1.} From this lemma   follow that
condition on fronts of the shock waves easy to get, considering
classical solutions of the hyperbolic equations as generalized
one. It is enough to equate  to zero the density,  corresponding
to independent singular generalized functions -  analogues of
simple, double and others layers, appearing under generalized
differentiation of the solutions.  The determination of such
conditions on base of the classical methods is more
labour-consuming procedure.

\textit{R e m a r k \,\,  2. } Eq.(\ref{(1.1)}) allows to consider
the generalized solutions with  derivatives breakup also on moving
surface $F(x,t) = 0$, which velocity of the motion can depend on
point of the front, then $v(x,t) = - F,_t /\left\| {grad\,F}
\right\|$ . On they the Adamar's  conditions (\ref{(1.4)})are
fulfilled with change $c$ to $v(x,t)$. Such solutions can be
generated by the right part of equation if the support of function
$G(x,t)$ enlarges in the course of time in $R^N $.

\vspace{3mm}

\textbf{2. Statement of nonstationary BVP.  Uniqueness of the
solution}. Let the solution of Eq.(\ref{(1.1)}) is determined in $
S^ - \subset R^N $, bounded by Lyapunov's surface $S$ ([2], p.
409), $t \ge 0$.\\
 \textit{Initial condition}: by $t = 0$
\begin{equation}\label{(2.1)}
u(x,0) = u_0 (x)\quad \textrm{for} \,\, x \in S^ - + S, \quad u,_t
(x,0) = \dot u_0 (x) \quad \textrm{ for} \,\, x \in S^-.
\end{equation}
We'll consider two BVP problems, consequently Dirihlet or Neumann
\textit{boundary conditions}:
\begin{equation}\label{(2.2)}
 u = u_S (x,t)\quad
\textrm{ for} \,\,x \in S,\;t \ge 0\quad\quad\textrm{(first BVP)};
\end{equation}
\begin{equation}\label{(2.3)}
\frac{{\partial u}}{{\partial n}} = p(x,t)\quad\quad \textrm{ for
}\,\, x \in S,\;t \ge 0\quad\quad \,\,\textrm{(second BVP)}.
\end{equation}
Here  we  denote as $n = \left( {n_1 ,...,n_N } \right)$ a unit
vector of external normal to $S$, $D^ - = S^ - \times R^ +$, $R^ +
= [0,\infty )$.

It is supposed that initial and boundary functions $u_0(x)$ and
$u_S(x,t)$ are continues,  $\dot{u}_0(x)$ and $p(x,t)$ piecewise
continues. For the first BVP the boundary and initial conditions
are coordinated:
\begin{equation}\label{(2.4)}
u_0 (x) = u_S (x,0)\quad \textrm{for} \,\, x \in S.
\end{equation}
On wave fronts, if they appear, Adamar's conditions (\ref{(1.4)})
are fulfilled.

Notice that shock waves always appear if there is not coordination
condition between initial velocities and  velocities at boundary:
\begin{equation}\label{(2.1*)}
\dot u_0 (x) = \dot u_S (x,0)\quad \textrm{for}\,\, x \in S,
\end{equation}
That is typically for physical problems. In this case at initial
moment of time at the boundary$D$  the shock wave front  is
formed, which spreads with the speed  $c$ in $R^N$. For building
continuously differentiated solutions this condition is necessary.
Here we don't suppose it.

It's supposed that initial conditions have been  given and the one
of the boundary conditions is known accordingly considered
problem.

We enter the functions $\;E = 0,5\left( {u,_\tau ^2 +
\sum\limits_{j = 1}^N {u,_j^2 } } \right)\;,\;L = 0,5\left(
{u,_\tau ^2  - \sum\limits_{j = 1}^N {u,_j^2 } } \right)$.

L e m m a \,\, 2.1. \textit{If $u(x,t)$ is a classical solution of
Eq. (\ref{(1.1)}) then }
\begin{equation}\label{(2.5)}
\left[ E \right]_{F_t }  =  - c^{ - 1} \left[ {\dot
u\frac{{\partial u}}{{\partial n}}} \right]_{F_t }, \quad \left[
{L(x,t)} \right]_{F_t }  =c ^{ - 2} \left( {\dot u^ -   +
c\,\frac{{\partial u^ -  }}{{\partial n}}} \right)\left[ {\;\dot
u} \right].
\end{equation}

P r o o f. On the strength of equality $[ab] = a^ + [b] + b^ -
[a]$ , with regard to (\ref{(1.1)}) and (\ref{(1.6*)}), we have
\[\begin{gathered}
 \left[ {cE + \dot u\frac{{\partial u}}{{\partial
n}}} \right] = \;\left[ {0,5\left( {c^{ - 1} \dot u^2  + cu,_j
u,_j } \right)\; + \dot u\frac{{\partial u}}{{\partial n}}}
\right] =\\ =  0,5c^{ - 1} \left( {\dot u^ +  \left[ {\dot u +
c\frac{{\partial u}}{{\partial n}}} \right] + \left( {\dot u^ - +
cu,_j^ -  n_j } \right)\left[ {\dot u} \right]} \right) + 0,5u,_j^
+  \left[ {\;cu,_j  + \dot un_j } \right]+\cr  + 0,5\left( {cu,_j^
- + \dot u^ -  n_j } \right)\left[ u,_j  \right] =  0,5c^{ - 1}
\left[ \dot u \right]\left( {\dot u^ {-}   + cu,_j^ -  n_j }
\right) + \\ + 0,5\left[ {u,_j } \right]\left( {cu,_j^ -   + \dot
u^ - n_j } \right) =   0,5cu,_j^ -  [ u,_j  + c^{ - 1} n_j \dot
{u} ] + 0,5c^{ - 1} \dot u^ -  [{cn_j u,_j  + \dot u} ] = 0
\end{gathered}
\]
(here $n$ is a unit normal to front $F_t$ in $R^N $). Thence  the
first formula of the lemma follows.

Hereinafter since $[a^2 ] = \left( {a^ + + a^ - } \right)[a]$ ,
and on the strength of (\ref{(1.4)}) and (\ref{(1.6)}), we get
second formula (\ref{(2.5)}):
\[
\displaylines{ \;\left[ L \right] = 0,5\left[ {u,_\tau ^2  -
\sum\limits_{j = 1}^N {u,_j^2 } } \right] = \;0,5\left( {u,_\tau ^
+   + u,_\tau ^ -  } \right)\left[ {u,_\tau ^{} } \right] -
0,5\left( {u,_j^ +   + u,_j^ -  } \right)\left[ {\;u,_j } \right]
=  \cr = 0,5\left( {u,_\tau ^ +   + u,_\tau ^ -  } \right)\left[
{u,_\tau ^{} } \right] + 0,5\left( {u,_j^ +   + u,_j^ -  }
\right)n_j \left[ {\;u,_\tau  } \right] =  \cr = 0,5\left\{
{\left( {u,_\tau ^ +   + n_j u,_j^ +  } \right) + \left( {u,_\tau
^ -   + n_j u,_j^ -  } \right)} \right\}\left[ {\;u,_\tau  }
\right] = c^{ - 2} \left( {\dot u^ -   + c\,\frac{{\partial u^ -
}}{{\partial n}}} \right)\left[ {\;\dot u} \right] .}
\]

\textit{R e m a r k}. If before the front of the wave $u \equiv
0$, that, with use of (\ref{(1.4*)}), we have $\left[ {L(x,t)}
\right]_{F_t } = 0$ i.e. in this case function $L$ continues.

T h e o r e m  2.1. \textit{If $u(x,t)$ is classical solution of
 BVP, then} $$ \int\limits_{S^ - } {(E(x,t)} - E(x,0))dV(x) = -
\int\limits_0^t {dt} \int\limits_{D^ -  } {G(x,t)u,_t dV(x)}  +
\int\limits_0^t {\int\limits_S { {\dot u_S (x,t)p(x,t)} dS(x)dt}
}.$$

P r o o f.  Multiplying Eq.(\ref{(1.1)}) on $u,_\tau $ in the
field of differentiability, after simple transformations we get:
\begin{equation}\label{(2.7)}
E,_\tau   - (u,_\tau  u,_j ),_j  =  - u,_\tau  G.
\end{equation}
But now we integrate Eq.(\ref{(2.7)}) on $D^ - $, with regard to
dividing  the area of integration $D^-$ on the parts, bounded
inside by wave fronts $F_k $.

Let us consider the  left part of this equation as divergency of
some vector in  $R^{N + 1} $, which  is continues and
differentiable in region between fronts. Then, using
Ostrogradskiy-Gauss theorem in $R^{N + 1} $, we get
\[
\displaylines{ \int\limits_{D^ -  } {E,_\tau  } dV(x)d\tau  -
\int\limits_{D^ -  } {(u,_\tau  u,_j ),_j dV(x)d\tau }  +
\int\limits_{D^ -  } {u,_\tau  G(x,\tau )dV(x)d\tau }  =  \cr =
\int\limits_{D^ -  } {u,_\tau  G(x,\tau )dV(x)d\tau }  +
\int\limits_{S^ -  } {(E(x,\tau )}  - E(x,0))dV(x) -
\int\limits_0^\tau  {\int\limits_S {\left( {u,_\tau  u,_j n_j }
\right)dS(x)d\tau } }  +  \cr + \sum\limits_{F_k }
{\int\limits_{F_k } {\left[ {E\nu _\tau   - u,_\tau  u,_j \nu _j }
\right]_{F_k } dF_k (x,\tau } } ) = 0,}
\]
here $dF_k (x,\tau )$ is a differential of the  area of surface in
corresponding point of  wave front $F_k$.   On the strength of
formulas (\ref{(1.3)}) and (\ref{(2.5)}), $\left[ {E\nu _\tau -
u,_\tau u,_j \nu _j } \right]_{F_k }  =  - c^{-1}{{\left\| \nu
\right\|_N }}\left[ {cE + \dot u\frac{{\partial u}}{{\partial n}}}
\right] = 0$. So the last integral is equal to zero. With regard
for indications for boundary function, thence the formula of the
theorem follows.

\textit{R e m a r k}.  First condition of formulas (\ref{(2.4)})
easy be get, considering Eq. (\ref{(2.7)}) in $D'(R^{N + 1} )$:\\
$ \hat E,_\tau   - (u,_\tau u,_j ),_j  =  - u,_\tau G + \left\{
{\left[ E \right]\nu _\tau   - \left[ {u,_\tau u,_j } \right]\nu
_j } \right\}\delta _F  =  - u,_\tau G - \left\| \nu \right\|_N
\left\{ {[ E ] + [ {u,_\tau u,_j } ]n_j }
\right\}\delta_{F_{t}}$.\\ Eq. (\ref{(2.5)}) is fulfilled in
$D'(R^{N + 1} )$ only if $ \left[ {E + c^{ - 1} \dot
u\frac{{\partial u}}{{\partial n}}} \right]_{F_t} = 0$

T h e o r e m 2.2. \textit{If classical solution of the first
(second) BVP exists, then it is single. }

P r o o f. On the strength of linearity of the problem, it is
enough to prove uniqueness of null solution. For it $G=0$, initial
and corresponding  boundary conditions are null. Then, as easy
see, from theorem 2.1 follows that $\int\limits_{S^ - } {E(x,t)}
dV(x) = 0$. Since $E$ is nonnegative, consequently $E \equiv 0
\Rightarrow u = const.$ From initial conditions follows $u \equiv
0.$

T e o r e m  2.3. \textit{If $u(x,t)$ is classical solution of the
BVP, then}
\[
\displaylines{ \int\limits_{D^ -  } {L(x,t)} dV(x)dt =
\int\limits_{D^ -  } {uG(x,t)dV(x)dt}  - \int\limits_0^t
{\int\limits_S {{u_S (x,t)p(x,t)} dS(x)dt} }  +  \cr + c^{ - 2}
\int\limits_{S^ -  } {\left( {u\dot u(x,t) - u_0 \dot u_0 (x)}
\right)dV(x)}  \cr}
\]
P r o o f. Multiplying Eq.(\ref{(1.1)}) to $u$, after simple
transformations we have:
\begin{equation}\label{(2.8)}
L(x,\tau ) + (uu,_j ),_j  - (uu,_\tau  ),_\tau   = G.
\end{equation}
Let integrate (\ref{(2.8)}) over $D^ - $ with regard for its
partition with wave fronts $F_k $. Similarly, as in theorem 2.1,
using Ostrogradskiy-Gauss theorem , we get
\[
\displaylines{ \int\limits_{D^ -  } {L\,} dV(x)d\tau  =
\int\limits_{D^ -  } {uGdV(x)d\tau }  + \int\limits_{D^ -  }
{\left( {(uu,_\tau  ),_\tau   - (uu,_j ),_j } \right)dV(x)d\tau }
=  \cr = \int\limits_{D^ -  } {uG\,dV(x)d\tau }  + c^{ - 1}
\int\limits_{S^ -  } {\left( {uu,_t  - u_0 \dot u_0 }
\right)dV(x)}  - \int\limits_0^\tau  {d\tau \int\limits_S {uu,_j
n_j dS(x)} }  +  \cr + \sum\limits_{F_k } {\int\limits_{F_k }
{\left[ {uu,_\tau  \nu _\tau   - uu,_j \nu _j } \right]_{F_k }
dF_k (x,\tau } } ).}
\]
On the strength of Adamar's conditions (\ref{(1.4)}),the last
summand is equal to zero. So, with regard for conditions on
boundary(\ref{(2.2)}),  (\ref{(2.3)}), we get the formula of the
theorem.

\textit{R e m a r k.} Theorem 2.1. is the law of energy
conservation. Similarly it's convenient  to use such law with
regard for conditions on front of the shock waves for proof of
uniqueness of BVP solutions in nonstationary problems of
mathematical physics .\vspace{3mm}

\textbf{3.  The dynamic analogue of  Green formula  in $D'(R^{N +
1} ).$ } For building the BVP solution we  consider it in $D'(R^{N
+ 1} )$. For this we  introduce the characteristic function of
semy\-cy\-lin\-der $D^ - $: $H_D^ - (x,t) \equiv H_S^ - \left( x
\right)H(t)$, where $H_S^ - \left( x \right)$is characteristic
function of set $S^ - $, which is equal to 0,5 on $S$, $H(t)$ is
Heaviside function equal to 0,5 by $t=0$. It's easy to show that
\begin{equation}\label{(3.1)}
\frac{{\partial H_D^ -  }}{{\partial x_j }} =  - n_j \delta _S
(x)H(t),\quad \frac{{\partial H_D^ -  }}{{\partial t}} =  - n_j
H_S^ -  (x)\delta (t).
\end{equation}
Hereinafter we will consider the generalized functions $$\hat
u(x,t) = u(x,t)H_D^ - (x,t),\quad\hat G(x,t) = G(x,t)H_S^ - \left(
x \right)H(t),$$ where $u(x,t)$ is the classical solution of BVP.

Let us define the action of wave operator on $\hat u(x,t)$. Since
$[u]_D = - u,$ executing generalized differentiation and using Eq.
(\ref{(1.1)}) in the field of differentiability, we get
\begin{equation}\label{(3.2)}
\begin{gathered}
\Box_c \hat u\left( {x,t} \right) =  - \frac{{\partial
u}}{{\partial n}}\delta _S \left( x \right)H\left( t \right) -
H\left( t \right)\left( {un_j \delta _S \left( x \right)}
\right),_j -\\ - c^{ - 2} H_S^ -  \left( x \right)u_0 \left( x
\right)\dot \delta (t) - c^{ - 2} H_S^ -  \left( x \right)\dot u_0
\left( x \right)\delta (t) + \hat G,
\end {gathered}
\end{equation}
where $\beta \left( {x,t} \right)\delta _S (x)H(t)$ is the simple
layer on  $\left\{ {S \times R^ + } \right\}$, $\delta (t)$ is
Dirak function , $\frac{{\partial u}}{{\partial n}} = u,_{i} n_i $
is a derivative along normal vector $n$ to $S$. Notice that
density of simple and double layers are here defined by use of
boundary conditions (a part from which, in depend on solved BVP,
are known) and by given initial conditions.

We can present the solution of the Eq. (\ref{(3.2)}) as the
convolution of  its right part  with the fundamental solution
$\hat U\left( {x,t} \right)$:
\begin{equation}\label{(3.3)}
\Box_c \hat U = \delta (x)\delta (t),
\end{equation}
\begin{equation}\label{(3.4)}
\hat U = 0\quad \textrm{for} \,\,t < 0 \,\,\textrm{and}\,\,
\left\| x \right\| > ct\quad(\textit{radiation conditions}).
\end{equation}
We shall name it \textit{Green function} of Eq. (\ref{(1.1)}).

The solution (\ref{(3.1)}) is represented in the manner of the
following convolution (*):
\begin{equation}\label{(3.5)}
\begin{gathered}
\hat{u}\left( {x,t} \right) = u\left( {x,t} \right)H_S^ - \left( x
\right)H\left( t \right) =  - \hat U\left( {x,t} \right)
* \frac{{\partial u}}{{\partial n}}\delta _S \left( x
\right)H\left( t \right) - \left( {\hat U * un_j \delta _S \left(
x \right)H\left( t \right)} \right),_j -\cr - c^{ - 2} \left(
{\hat U\mathop  * \limits_x H_S^ -  \left( x \right)u_0 \left( x
\right)} \right),_t  - c^{ - 2} \hat U\mathop  * \limits_x H_S^ -
\left( x \right)\dot u_0 \left( x \right) + \hat G * \hat U,\\
\end{gathered}
\end{equation}
where symbol "$\mathop * \limits_x $" means that convolution  on
$x$ is taken only  over $x$. Moreover the solution (\ref{(3.5)})
is single solution in the class of functions, allowing the
convolution with $U(x,t)$.

If initial data and $G $ are null then
\begin{equation}\label{(3.5*)}
\hat u\left( {x,t} \right) =  - \hat U\left( {x,t} \right) *
\frac{{\partial u}}{{\partial n}}\delta _S \left( x \right)H\left(
t \right) - \hat U,_j  * un_j \delta _S \left( x \right)H\left( t
\right).
\end{equation}

This formula expresses the BVP solution  through boundary values
of  function $u$ and its normal derivative and similar to Green
formula  for solutions of  Laplace equation [2]. However, on the
strength of particularities of the fundamental solutions of
hyperbolic equations on  wave fronts, which type depends on space
dimensionality , its integral representation gives the dispersing
integrals (in second summand).

For building its regular integral presentation we introduce the
antiderivatives on $t$ functions:
\begin{equation}\label{(3.6)}
\hat W\left( {x,t} \right) = \hat U\left( {x,t} \right) * \delta
\left( x \right)H\left( t \right) = \hat U\left( {x,t}
\right)\mathop  * \limits_t H\left( t \right)\,\, \Rightarrow \,\,
\partial _t \hat W\left( {x,t} \right) = \hat U\left( {x,t} \right)
;
\end{equation}
\[
\hat H\left( {x,m,t} \right) = \frac{{\partial \hat W\left( {x,t}
\right)}}{{\partial x_i }}m_i  = \frac{{\partial \hat W\left(
{x,t} \right)}}{{\partial m}}.
\]
Easy to see that they also are  solutions of Eq.(\ref{(1.1)}) when
$G = H(t)\delta (x)$ and $G = H(t)\frac{{\partial \delta \left(
{x,t} \right)}}{{\partial m}}$ accordingly.

T h e o r e m  3.1. \textit{In $D'(R^{N + 1} )$ BVP solution
satisfies to the equation}:
\begin{equation}\label{(3.7)}
\hat u\left( {x,t} \right) =  - \hat U\left( {x,t} \right) *
\frac{{\partial u}}{{\partial n}}\delta _S \left( x \right)H\left(
t \right) - \left( {\hat W,_j  * \dot un_j \delta _S \left( x
\right)H\left( t \right)} \right) -
\end{equation}
\[
- {\hat W,_j \mathop * \limits_x u_0 \left( x \right)n_j \left( x
\right)\delta _S \left( x \right) - c^{ - 2} \hat U\mathop  *
\limits_x H_S^ -  \left( x \right)\dot u_0 \left( x \right) - c^{
- 2} \left( {\hat U\mathop * \limits_x H_S^ - \left( x \right)u_0
\left( x \right)} \right),_t + \hat G * \hat U} yu
\]

P r o o f. It's easy to show, using determination of derivative of
generalized functions and continuity $u$, that
\[
\left( {un_j \delta _S \left( x \right)H\left( t \right)}
\right),_t  = \dot u(x,t)n_j (x)\delta _S \left( x \right)H\left(
t \right) + u(x,0)n_j (x)\delta _S \left( x \right)\delta \left( t
\right).
\]
Using  this relation, last formula(\ref{(3.6)}) and  properties of
convolution differentiation , we have
\[
\left( {\hat U * un_j \delta _S \left( x \right)H\left( t \right)}
 \right),_j   = \hat W,_j  * \dot u(x,t)n_j (x)\delta _S \left( x \right)H\left( t \right)
  + \hat W,_j \mathop * \limits_x u_0 (x)n_j
(x)\delta _S \left( x \right)
\]
Substituting these correlations in formula (\ref{(3.5)}), we get
the formula of the theorem.

From theorem follow that solution of the problem is completely
defined by boundary value of the normal derivatives of $u\left(
{x,t} \right)$ and its velocities $\dot u $. In analogy with
presentation of the solutions of the Laplace equation, these
formulas possible to name \textit{dynamic analogue of Green
formula}.

The formula (\ref{(3.7)}) possesses the advantage in comparison
with (\ref{(3.5)}), since it allows immediately to go to its
integral record without regularization of integral function on
front, as it was  earlier offered in paper [3]. For $x \in S$
formula (\ref{(3.7)}) gives, what we show hereinafter, boundary
integral equation for BVP solution. If  one of the boundary
function is known, then after solving BIE on border, we can find
second boundary function. Whereupon formula (\ref{(3.7)}) gives
the solution $u(x,t)$ in $D^-$.\vspace{3mm}

\textbf{4. Dynamic analogue of Gauss formula.} We introduce the
functions
\[
U\left( {x,y,t} \right) = \hat U\left( {x - y,t} \right) ,\,\,
W\left( {x,y,t} \right) = \hat W\left( {x - y,t} \right) ,\,\,
H\left( {x,y,m,t} \right) = \hat H\left( {x - y,m,t} \right),
\]
which, on the strength of properties of  symmetries of the wave
operator and $\delta$-functions, satisfy to following symmetries
correlations:
\begin{equation}\label{(4.1)}
\begin{gathered}
U\left( {x,y,t} \right) = U\left( {y,x,t} \right) ,\,\, W\left(
{x,y,t} \right) = W\left( {y,x,t} \right) ,\quad \frac{{\partial
W}}{{\partial x_j }} =  - \frac{{\partial W}}{{\partial y_j }},
\cr H\left( {x,y,m,t} \right) =  - H\left( {y,x,m,t} \right) = -
H\left( {x,y, - m,t} \right).
\end{gathered}
\end{equation}

L e m m a 4.1. \textit{In $D'(R^{N + 1} )$ the dynamic analogue of
Gauss formula  has the form}:
\begin{equation}\label{(4.2)}
- \hat W,_i \mathop * \limits_x n_i \left( x \right)\delta \left(
x \right) - c^{ - 2} \left( {\hat U\mathop  * \limits_x H_S^ -
\left( x \right)} \right),_t  = H_D^ -  \left( {x,t} \right)
\end{equation}
P r o o f. If in both parts of Eq. (\ref{(3.3)}) for $U$ the
convolution with $H_D^-(x,t)$ are taken, after using a property of
convolution differentiation  and formula (\ref{(3.1)}) we get:
\[
- U,_j *n_j (x)\delta _S \left( x \right)H(t) - c^{ - 2} \left(
{U,_t *H_S^ -  \left( x \right)\delta (t)} \right) = H_D^ -  (x,t)
\]
With regard for (\ref{(3.6)}), flipping differentiation on $t$ in
the first summand and executing convolution on $t$ in second one,
we get the formula of this lemma. Integral record of this formula
depends on dimensionality of the equation.

The formula (\ref{(4.2)}) is an analogue of the known formula of
Gauss for  double layer potential (p.406, [2]), which gives the
integral record of characteristic function of a set with use the
fundamental solution of  Laplace equation. The Gauss formula is
often used for building the boundary integral equation of boundary
value problems for elliptical equations and systems. Similarly
it's possible to use the dynamic analogue of Gauss formula  for
constructing BIE in hyperbolic case. But here the dynamic analogue
of  Green formula  is used for building BIE of setting problems.

Hereinafter we shall give integral representation of formulas of
the theorem 3.1. and lemma 4.1 for space of dimensionality
N=1,2,3, that is most typical for problems of mathematical
physics.\vspace{3mm}

\textbf{5. BIE for plane BVP}.  When $N=2$  we have a plane
problem. At first let's consider the  problem with null initial
conditions.

We denote as $dS(y)$ the arc length differential on $S$ in point
$y$, $S_t (x) = \left\{ {y \in S:\, r  < ct} \right\}$, $S_t^ -
(x) = \{y \in S^ - : r < ct \},\, r = \left\| {x - y}
\right\|,\,dV(y) = dy_1 dy_2  $.

T h e o r e m \,\, 5.1. \textit{By $N=2$ function $\hat u\left(
{x,t} \right)$ with null initial conditions $\left( {u_0 =
0,\;\dot u_0 = 0} \right)$, has  following integral
representation}:
\begin{equation}\label{(5.1)}
\hat u = \frac{c}{{2\pi }}\int\limits_0^t {d\tau }
\int\limits_{S_\tau  ( x )} \left( {\frac{\partial u(y,t -
\tau)}{\partial n( y )} + \frac{\tau }{r}\frac{\partial
r}{\partial n (y)}\dot u(y,t - \tau ) }\right )\frac{dS(y)}{\sqrt{
c^2 \tau ^2 - r^2}} .
\end{equation}
\textit{Or, if to change the integration order, it has the form }:
\[
\begin{gathered}
\hat u\left( {x,t} \right) = \frac{c}{{2\pi }}\int\limits_{S_t
\left( x \right)} {dS\left( y \right)} \int\limits_{{r
\mathord{\left/ {\vphantom {r c}} \right.
\kern-\nulldelimiterspace} c}}^t {\frac{\partial u\left( {y,t -
\tau } \right)}{\partial n\left( y \right)}\frac{d\tau
}{\sqrt{\left( {c^2 \tau ^2  - r^2 } \right)}}}+ \\ +
\frac{c}{{2\pi }}\int\limits_{S_t \left( x \right)}
{\frac{1}{r}\frac{{\partial r}}{{\partial n\left( y
\right)}}dS\left( y \right)} \int\limits_{{r \mathord{\left/
{\vphantom {r c}} \right. \kern-\nulldelimiterspace} c}}^t
{\frac{{\tau \dot u\left( {y,t - \tau } \right)d\tau
}}{{\sqrt{\left( {c^2 \tau ^2  - r^2 } \right)}} }}.
 \end{gathered}
\]
\textit{For $x\in S$ second integral on the right part is singular
and taken in value principle sense.}

P r o o f. For $N=2$ Green function has a following type ([2], p.
206)
\begin{equation}\label{(5.2)}
\hat U\left( {x,t} \right) =  - \frac{{cH\left( {ct - R}
\right)}}{{2\pi {\sqrt{\left( {c^2 \tau ^2  - r^2 } \right)}} }}
,\quad R = \sqrt {x_1^2  + x_2^2 }.
\end{equation}
Calculating from formula (\ref{(3.6)}) we determine
\begin{equation}\label{(5.3)}
\hat W\left( {x,t} \right) =  - \frac{{H\left( {ct - R}
\right)}}{{2\pi }}\ln \left( {\frac{{ct - \sqrt {c^2 t^2  - R^2 }
}}{R}} \right) ,\quad \hat H\left( {x,m,t} \right) =  -
\frac{{\,\,ct\,H\left( {ct - R} \right)}}{{2\pi \,\sqrt {c^2 t^2 -
R^2 } }}\frac{{x_j m_j }}{{R^2 }} .
\end{equation}

If to write the convolutions (\ref{(3.7)}) in integral type with
regard for these relations, that we get the formula of the
theorem. Notice that for $ x \notin S$ integrals on the right are
regular and therefore for such $x$ on the right and on the left
the regular functions stand . We shall prove formula (\ref{(5.1)})
also for $x \in S$ with regard for determinations $H_S^ - (x)$.

We denote $\varepsilon$-vicinity of the point $x$ ($\varepsilon
\ll ct, t>0)$ through $ \textrm{Ш}_\varepsilon (x) = \left\{
{y:\;r < \varepsilon } \right\},$ $ S_\varepsilon ^ - \left( x
\right) = S^ - - \textrm{Ш}_\varepsilon \left( x \right),
S_\varepsilon ^ + \left( x \right) = S^ + - \textrm{Ш}_\varepsilon
\left( x \right),\,O _\varepsilon \left( x \right) = \left\{ {y
\in S:\;r \le \varepsilon } \right\}, \,\,S_\varepsilon   = S -
O_\varepsilon,$ $\textrm{Ш}_\varepsilon ^ - = S^ - \cap
\textrm{Ш}_\varepsilon,\,\,\textrm{Ш} _\varepsilon ^ + = S^ + \cap
\textrm{Ш}_\varepsilon,\,\,\textrm{Г} _\varepsilon ^ \pm \left( x
\right) = \left\{ {y \in S^ \pm  :\;r = \varepsilon } \right\},
\quad r,_j = \frac{{\partial r}}{{\partial y_j }}.$

Let $x^*\in S$. We transform the sidebar $S$ in vicinity of the
point $x^*$, avoiding it on $\varepsilon $-semicircle in $S^{-}\,
(\varepsilon <<ct, t>0)$. From  dynamic analogue of the Green
formula  for sidebar $S_\varepsilon +\textrm{Г} _\varepsilon ^ - $
we have in point $x=x^*$
\[
0 = \frac{c}{\pi }\int\limits_{S_\varepsilon  \left( x^* \right) +
\textrm{Г}_\varepsilon ^ - \left( x^* \right)} {H\left( {ct - r}
\right)dS\left( y \right)\int\limits_{{r \mathord{\left/
{\vphantom {r c}} \right. \kern-\nulldelimiterspace} c}}^t
{\frac{{\partial u\left( {y,t - \tau } \right)}}{{\partial n\left(
y \right)}}\frac{{d\tau }}{{\sqrt{c^2 \tau ^2  - r^2}} }}} +
\]
\[
+ \frac{c}{{2\pi }}\int\limits_{S_\varepsilon  \left( x^* \right)}
{H\left( {ct - r} \right)\frac{1}{r}\frac{{\partial r}}{{\partial
n\left( y \right)}}dS\left( y \right)\int\limits_{{r
\mathord{\left/ {\vphantom {r c}} \right.
\kern-\nulldelimiterspace} c}}^t {\frac{{\tau \dot u\left( {y,t -
\tau } \right)d\tau }}{{\sqrt{c^2 \tau ^2  - r^2}} }}} +
\]
$$
+ \frac{c}{{2\pi }}\int\limits_{\textrm{Г}_\varepsilon ^ - \left(
x^* \right)} {H\left( {ct - r} \right)\frac{1}{r}\frac{{\partial
r}}{{\partial n\left( y \right)}}dS\left( y \right)\int\limits_{{r
\mathord{\left/ {\vphantom {r c}} \right.
\kern-\nulldelimiterspace} c}}^t {\frac{{\tau \dot u\left( {y,t -
\tau } \right)d\tau }}{\sqrt{c^2 \tau ^2  - r^2} }}}.
$$
For $\varepsilon \to {\rm + 0}$ first integral, on the strength of
weak singularity  of the integral functions, strives to integral
on $S$, second one tends to integral in the value principle sense,
which also exists, since integral function has a singularity of
the type $r^{ - 1} $ and contains function $\frac{{\partial
r}}{{\partial n(y)}} = n_j (y)\frac{{\left( {y_j - x_j }
\right)}}{r}$, which for $y \to x$ asymptotically equivalent to
$\frac{{\partial r}}{{\partial n(x)}} = n_j (x)\frac{{\left( {y_j
- x_j } \right)}}{r}$,  antisymmetric relative to point $x$.

Let consider the last summand,  denoted $J_\varepsilon \left( x
\right)$. On $\textrm{Г}_\varepsilon ^ - $ there are $ r =
\varepsilon ,\, \frac{{\partial r}}{{\partial n(y)}} =   - 1,\,
dS\left( y \right) = \varepsilon \,d\theta$ , where $\theta$ is an
polar corner with top in point $x$; $\theta _1 $ and $\theta _2 $
are the corners of endpoints  $\textrm{Г}_\varepsilon ^ - $,
numbered in order at pass-by of the sidebar
$\textrm{Г}_\varepsilon ^{-} $ in positive direction. With regard
for these relations
$$
\begin{gathered}
J_\varepsilon  \left( x \right) = \frac{c}{{2\pi
}}\int\limits_{\theta _1 }^{\theta _2 } {\frac{\varepsilon
}{\varepsilon }d\theta } \int\limits_{{\varepsilon \mathord{\left/
{\vphantom {\varepsilon  c}} \right. \kern-\nulldelimiterspace}
c}}^t {\frac{\tau \dot u\left( {y,t - \tau } \right)d\tau
}{\sqrt{c^2 \tau ^2  - r^2}}}= \frac{{\left( {\theta _2 - \theta
_1 } \right)c}}{{2\pi }}\int\limits_{{\varepsilon \mathord{\left/
{\vphantom {\varepsilon  c}} \right. \kern-\nulldelimiterspace}
c}}^t {\frac{{\tau \dot u\left({y,t - \tau } \right)d\tau
}}{\sqrt{c^2 \tau ^2  - r^2}} }, \cr \mathop {\lim
}\limits_{\varepsilon  \to 0} \left( {\theta _2  - \theta _1 }
\right) =  - \pi ,\quad \mathop {\lim }\limits_{\varepsilon \to 0}
J_\varepsilon  \left( x \right) =  - 0,5\int\limits_0^t {\dot
u\left( {y,t - \tau } \right)d\tau  =  - 0,5u\left( {x,t}
\right)}. \end{gathered} $$

Passing this summand in the left part, with regard for
determinations: $H_D^ - (x,t) = 0,5$ for $x \in S$, we get  the
formula of theorem on the border also. Since on the left and on
the right in Eq. (\ref{(5.1)}) the regular generalized functions
stand, on the strength of Dyubua-Reymon lemma ([2], p. 97) this
equality, which equitable in generalized functions class , also
equitable in classical sense. This theorem has been proved.

By solving  the BVP the Eq. (\ref{(5.1)}) gives the boundary
integral equation for determination of  unknown boundary function:
\[
\begin{gathered}
\pi  u_S \left( {x,t} \right) = \\=\int\limits_{S_t \left( x
\right)} {dS\left( y \right)} \int\limits_{{r \mathord{\left/
{\vphantom {r c}} \right. \kern-\nulldelimiterspace} c}}^t
{\frac{{\partial u\left( {y,t - \tau } \right)}}{{\partial n\left(
y \right)}}\frac{{d\tau }}{\sqrt{ \tau ^2  - r^2/c^2}} } +
V.P.\int\limits_{S_t \left( x \right)} {\frac{1}{r}\frac{{\partial
r}}{{\partial n\left( y \right)}}dS\left( y \right)}
\int\limits_{{r \mathord{\left/ {\vphantom {r c}} \right.
\kern-\nulldelimiterspace} c}}^t {\frac{{\tau \dot u_S \left( {y,t
- \tau } \right)d\tau }}{\sqrt{\tau ^2 - r^2/c^2}} }.
\end{gathered}
\]
In the case of the first BVP  the left part of this equation and
the second integral on the right are known, they are defined with
use of boundary condition, but the first integral contains the
kernel with weak singularity on front of the Green function.
Solving it we define the normal derivative of $u$ on boundary
whereupon formula (\ref{(5.1)}) allows to compute the solution in
any point.

In the case of the second BVP we have \textit{singular} BIE for
determination of unknown boundary values of  $u$. After its
solving we can define $u(x,t)$ on boundary whereupon formula
(\ref{(5.1)}) defines the solution fully.

In the case of nonzero initial conditions the following theorem
gives the solution of problems.

T h e o r e m  5.2. \textit{By $N=2$ BVP solution has the
following integral representation}:
\[
\displaylines{ 2\pi \hat u = \int\limits_0^t {d\tau }
\int\limits_{S_\tau  \left( x \right)} {\left( {\frac{{\partial
u\left( {y,t - \tau } \right)}}{{\partial n\left( y \right)}} +
\frac{1}{r}\frac{{\partial r}}{{\partial n\left( y \right)}}\tau
\dot u\left( {y,t - \tau } \right)} \right)} \frac{{dS\left( y
\right)}}{{\sqrt {\tau ^2  - \left( {r/c} \right)^2 } }} +  \cr +
\frac{\partial }{{\partial t}}\int\limits_{S_t^ -  (x)}
{\frac{{u_0 (y)dV\left( y \right)}}{{c\sqrt {c^2 t^2  - r^2 } }}}
+\int\limits_{S_t^ -  (x)} {\frac{{\dot u_0 (y)dV\left( y
\right)}}{{c\sqrt {c^2 t^2  - r^2 } }}}  + c\int\limits_0^t {d\tau
} \int\limits_{S_\tau^-  \left( x \right)} {\frac{{G(y,t - \tau
)dV\left( y \right)}}{{\sqrt {c^2 \tau ^2  - r^2 } }}}  +  \cr -
\int\limits_{S_t (x)} {u_0 (y)\frac{{ct}}{r}} \frac{{\partial
r}}{{\partial n\left( y \right)}}\frac{{dS\left( y
\right)}}{{\sqrt {c^2 \tau ^2  - r^2 } }} }
\]
Proof  follows from theorem 4.1. and 5.2. if to write the
convolution with initial data  in integral type. Here the
integrals from the second before the fourth ones comply with
Poisson formula for Caushy problem . The last additional summand
with initial data is conditioned presence of the boundary.

Using relations (\ref{(5.2)}), dynamic analogue of Gauss formula
possible to be write in integral type.

L e m m a 5.1. \textit{By $N=2$ dynamic analogue of  Gauss formula
has the following type}:
\[
V.P.\int\limits_{S_\tau  \left( x \right)} {\frac{1}{{\sqrt {1 -
\left( {r/ct} \right)^2 } }}\frac{\partial }{{\partial n\left( y
\right)}}\left( {\ln \frac{1}{r}} \right)dS\left( y \right) +
\frac{\partial }{{c\partial t}}\int\limits_{S_t^ -  \left( x
\right)} {\frac{{dV\left( y \right)}}{{\sqrt {c^2 t^2  - r^2 } }}}
= 2\pi H_S^ - ( x )H( t ),}
\]
\textit{where integral in the value principle sense is taken  for
boundary points }.

P r o o f. Formula of the lemma 4.1 , with regard for
(\ref{(5.2)}),(\ref{(5.3)}), possible to write so
\begin{equation}\label{(5.4)}
- \int\limits_{S_t \left( x \right)} {\frac{{ct}}{{\sqrt {c^2 t^2
- r^2 } }}\frac{1}{r}\frac{\partial r}{\partial n\left( y
\right)}dS(y) + \frac{\partial }{{c\partial t}}\int\limits_{S_t^ -
(x)} {\frac{{\,dV(y)}}{{\sqrt {c^2 t^2 - r^2 } }}}}= 2\pi H_S^ -
(x)H(t).
\end{equation}
From this formula by elementary transformations we get the formula
of the lemma.

We shall show that this equality, equitable in the field of
regularity , is saved also for $x \in S$ if the singular integral
on the left, which contains the strong singularity on $r$, is
taken in the value principle sense.

Similarly (\ref{(5.4)}) for area without and with
 $\varepsilon$ -vicinity of the point $x$ we get
\begin{equation}\label{(5.5)}
\begin{gathered}
- \int\limits_{S_\varepsilon + \textrm{Г}_\varepsilon ^ - }
{\frac{ct\;H\left( {ct - r} \right)}{\sqrt{c^2 \tau ^2  - r^2}}
\frac{1}{r}\frac{{\partial r}}{{\partial n\left( y
\right)}}dS\left( y \right) + \frac{\partial }{{c\partial
t}}\int\limits_{S_\varepsilon ^ - \left( x \right)}
{\frac{{H\left( {ct - r} \right)dV\left( y \right)}}{\sqrt{c^2
\tau ^2  - r^2}} }} = 0,\\
- \int\limits_{S_\varepsilon + \textrm{Г}_\varepsilon ^ + }
{\frac{{ct\;H\left( {ct - r} \right)}}{\sqrt{c^2 \tau ^2  - r^2}}
\frac{1}{r}\frac{{\partial r}}{{\partial n\left( y
\right)}}dS\left( y \right) + \frac{\partial }{{c\partial
t}}\int\limits_{S_\varepsilon ^ - \left( x \right) +
\textrm{Ш}_\varepsilon } {\frac{{H\left( {ct - r} \right)dv\left(
y \right)}}{\sqrt{c^2 \tau ^2  - r^2}} }}  = 2\pi.\\
\end{gathered}
\end{equation}
Under $\varepsilon < ct$ integrals on
$\textrm{Г}_{\varepsilon}^{\pm},\,\textrm{Ш}_\varepsilon$ are easy
calculated by transition to polar coordinate system. On
$\textrm{Г}_\varepsilon ^ \pm $ $\frac{{\partial r}}{{\partial
n(y)}} = \pm1$ consequently to the sign. Mark that
$$
\begin{gathered}
I_\varepsilon  \left( x,t \right) =
\int\limits_{\textrm{Ш}_\varepsilon \left( x \right)}
{\frac{{H\left( {ct - r} \right)dV\left( y \right)}}{\sqrt{c^2
\tau ^2  - r^2}} }  = \int\limits_0^{2\pi } {d\theta }
\int\limits_0^\varepsilon  {\frac{{rdr}}{\sqrt{c^2 \tau ^2  -
r^2}} }  = 2\pi (ct-\sqrt{c^2 \tau ^2  - r^2} ) \quad \Rightarrow
\\\mathop {\lim }\limits_{\varepsilon  \to 0} \frac{{\partial
I_\varepsilon  }}{{\partial t}} = 2\pi c\mathop {\lim
}\limits_{\varepsilon  \to 0} \left( {1 - \frac{ct}{\sqrt{c^2 \tau
^2 - r^2}}} \right) = 0\,\, \textrm {for} \,\,\forall t > 0.\\
\end{gathered}
$$
If both equalities in (\ref{(5.5)}) to add and divide on 2, with
regard for properties of the symmetries of integral function
(\ref{(4.1)}) and going to limit on $\varepsilon \to 0$ , then
$$
- V.P.\int\limits_{S_t \left( x \right)} {\frac{{ct\;}}{\sqrt{c^2
\tau ^2  - r^2}} }\frac{1}{r}\frac{{\partial r}}{{\partial n\left(
y \right)}}dS\left( y \right) + \frac{\partial }{{c\partial t}}
\int\limits_{S_t^ -  \left( x \right)} {\frac{{dV\left( y
\right)}}{\sqrt{c^2 \tau ^2  - r^2}}}  = \pi.
$$
So, with regard for determinations $H_S^ - (x)$, formula
(\ref{(5.4)}) is equitable for any $x$.\vspace{3mm}

\textbf{6. BIE of  BVP for $\textbf{N=3.}$}

T h e o r e m   6.1. \textit{For  $N=3$ BVP solution is
represented in the form}:
\[
\begin{gathered}
4\pi \hat{u}\left( {x,t} \right) = \int\limits_{S_t \left( x
\right)} {\left\{ {\frac{1} {r}\frac{{\partial u\left( {y,t - {r
\mathord{\left/ {\vphantom {r c}} \right.
\kern-\nulldelimiterspace} c}} \right)}} {{\partial n\left( y
\right)}} + \frac{{\dot u\left( {y,t - {r \mathord{\left/
{\vphantom {r c}} \right. \kern-\nulldelimiterspace} c}} \right)}}
{c}\frac{{\partial \,\ln \,r}} {{\partial n\left( y \right)}}}
\right\}} dS\left( y \right)-
\\
- V.P.\int\limits_{S_t \left( x \right)} {u\left( {y,t - {r
\mathord{\left/ {\vphantom {r c}} \right.
\kern-\nulldelimiterspace} c}} \right)\frac{\partial } {{\partial
n\left( y \right)}}\frac{1} {r}} dS\left( y \right) + c^{ - 1}
\frac{\partial } {{\partial t}}\int\limits_{S_t \left( x \right)}
{u_0 \left( y \right)\frac{{\partial \,\ln \,r}} {{\partial
n\left( y \right)}}dS\left( y \right)}  + \\+ {\int\limits_{r =
ct} {\frac{{\dot u_0 \left( y \right)}} {{c^2 t}}H_S^ -  \left( y
\right)dS\left( y \right)}  +   \frac{\partial } {{\partial
t}}\int\limits_{r = ct} {\frac{{u_0 \left( y \right)}} {{c^2
t}}H_S^ -  \left( y \right)dS\left( y \right)} }  -
\int\limits_{S_t^ -  (x)} {\frac{{G(x,t - r/c)}} {r}} dV(y).
\end{gathered}
\]

P r o o f . For $N=3$ Green function $\hat U(x,t)$ is the double
layer on the cone $K_t = \left\{ {\left( {x,t} \right):\;\left\| x
\right\| = ct} \right\}$:
\begin{equation}\label{(6.1)}
\hat U\left( {x,t} \right) =  - \frac{{\delta \left( {t - R/c}
\right)}} {{4\pi R}} , \quad R = \sqrt {x_1^2  + x_2^2  + x_3^2 }.
\end{equation}
For any $\varphi \left( {x,t} \right) \in D\left( {R_4 } \right)$
it defines the linear functional:
\[
\left( {U,\varphi } \right) =  - \frac{1} {{4\pi
}}\int\limits_{R^3 } {\frac{{\varphi \left( {x,\left\| x
\right\|/c} \right)}} {{\left\| x \right\|}}dV\left( x \right)}
\]
After calculation on formula (\ref{(3.6)}) we find
\begin{equation}\label{(6.2)}
\hat W\left( {x,t} \right) =  - \frac{{H\left( {ct - R} \right)}}
{{4\pi R}} ,\,\, \hat H\left( {x,m,t} \right) = \frac{1} {{4\pi
}}\frac{{x_j m_j }} {{R^2 }}\left( {c^{ - 1} \delta \left( {t -
R/c} \right) + \frac{{H\left( {ct - R} \right)}} {R}} \right).
\end{equation}

For building  the integral representation of dynamic analogue of
Green and Gauss formulas we use the following equalities, which be
get using the determination  of  generalized function convolution
:
\[
\alpha \left( x \right)\delta \left( {t - R/c} \right) * f\left(
{x,t} \right)H\left( t \right) = H\left( t \right)\int\limits_{S_t
\left( x \right)} {\alpha \left( {x - y} \right)f\left( {y,t - {r
\mathord{\left/ {\vphantom {r c}} \right.
\kern-\nulldelimiterspace} c}} \right)dS\left( y \right)}
\]
\[
\beta \left( {x,t} \right)\delta _S \left( x \right)H\left( t
\right) * f\left( {x,t} \right)H\left( t \right) = H\left( t
\right)\int\limits_0^t {d\tau } \int\limits_{S_t (x)} {\beta
\left( {y,t - \tau } \right)f\left( {x - y,\tau } \right)dS\left(
y \right)}
\]
\[
\alpha \left( x \right)\delta \left( {t - R/c} \right) * \beta
\left( {x,t} \right)\delta _S \left( x \right)H\left( t \right) =
H\left( t \right)\int\limits_{S_t \left( x \right)} {\alpha \left(
{x - y} \right)\beta \left( {y,t - {r \mathord{\left/ {\vphantom
{r c}} \right.\,\, \kern-\nulldelimiterspace} c}} \right)dS\left(
y \right)},
\]
\[
\alpha \left( x \right)\delta \left( {t - R/c} \right)\mathop  *
\limits_x \gamma \left( x \right)\delta _S \left( x \right) =
\frac{\partial } {{\partial t}}\int\limits_{S_t \left( x \right)}
{\alpha \left( {x - y} \right)\gamma \left( y \right)dS\left( y
\right)}
\]
\[
\alpha \left( x \right)\delta \left( {t - R/c} \right)\mathop  *
\limits_x H_S^ -  \left( x \right) = c^{ - 1} H\left( t
\right)\int\limits_{r = ct} {\alpha \left( {x - y} \right)H_S^ -
\left( y \right)dS\left( y \right)}.
\]
We compute the convolution in Eq.(\ref{(3.7)}) with regard for
these correlations. First summand in this equation is
\[
 - \hat U\left( {x,t} \right) * \frac{{\partial u}}
{{\partial n}}\delta _S \left( x \right)H\left( t \right)  =
\frac{{H\left( t \right)}} {{4\pi }}\int\limits_{S_t \left( x
\right)} {\frac{1} {r}\frac{{\partial u\left( {y,t - {r
\mathord{\left/ {\vphantom {r c}} \right.
\kern-\nulldelimiterspace} c}} \right)}} {{\partial n(y)}}dS\left(
y \right)}.
\]
The second summand ($ - 4\pi \,\hat W,_j * \dot u(x,t)n_j
(x)\delta _S \left( x \right)H\left( t \right)$) is an amount of
two convolutions:
\[
\begin{gathered}
- \frac{1} {{4\pi c}}\frac{{x_j }} {{R^2 }}\delta \left( {t -
\left\| x \right\|/c} \right) * \dot u(x,t)n_j (x)\delta _S \left(
x \right)H\left( t \right) = c^{ - 1} H(t)\int\limits_{S_t (x)}
{\dot u(y,t - r/c)} \frac{{\partial \ln r}}
{{\partial n(y)}}dS(y),\\
 - \frac{{x_j }} {{R^3 }}H\left( {ct -
\left\| x \right\|} \right)
* \dot u(x,t)n_j (x)\delta _S \left( x \right)H\left( t \right)=\\
= H(t)\int\limits_{S_t (x)} {u^0 (y)\frac{{\partial r^{ - 1} }}
{{\partial n(y)}}dS(y) - } H(t)\int\limits_{S_t (x)}
{u(y,t)\frac{{\partial r^{ - 1} }}
{{\partial n(y)}}dS(y)}  \\
\end{gathered}
\]
Third summand
\[
\begin{gathered}
- 4\pi \,\hat W,_j \mathop * \limits_x u_0 \left( x \right)n_j
\left( x \right)\delta _S \left( x \right) = \frac{\partial }
{{c\partial t}}\int\limits_{S_t \left( x \right)} {\frac{{\partial
\ln r}} {{\partial n\left( y \right)}}u_0 \left( y \right)dS\left(
y \right)}  - \int\limits_{S_t \left( x \right)} {\frac{{\partial
r^{ - 1} }}
{{\partial n\left( y \right)}}u_0 \left( y \right)dS\left( y \right)}  \\
\end{gathered}
\]
The last three summand give the known Kirhoff formula  for
solution of the Caushy problem [ 1,2] with initial conditions.
Summing these convolutions, we get the formula of the theorem.

In the formula of the theorem in the part, depending on initial
data, there is one summand conditioned by presence of the
boundary$S$. It disappears for $t > t^* (x),\, \,t^* (x) = \mathop
{\max }\limits_{y \in S} \frac{{\left\| {x - y} \right\|}}{c}$
since $S_t \left( x \right) = S$ and integral does not depend on $
t $.

In formula (\ref{(3.1)}) for $x \notin S,\, t>0$ all integrals
exist. Its proof  for $x \in S$ like to plane case. Herewith the
strong singularity has the second summand on the right. In this
case
\[
\mathop {\lim }\limits_{\varepsilon  \to 0}
\int\limits_{\textrm{Г}_\varepsilon ^ - } {u\left( {y,t - \frac{r}
{c}} \right)} \frac{\partial } {{\partial n\left( y
\right)}}\frac{1} {r}dS\left( y \right) = -\mathop {\lim
}\limits_{\varepsilon  \to 0} \frac{1} {{\varepsilon ^2
}}\int\limits_{\textrm{Г}_\varepsilon ^ - } {u\left( {y,t -
\frac{r} {c}} \right)} dS\left( y \right) =  - 2\pi u\left( {x,t}
\right),
\]
It's clear that for $x \in S$ formula (\ref{(6.1)}) saves kind if
corresponding singular integral to calculate in the value
principle sense and take into account the value $H_S^ - \left( x
\right)$ on $S$.

For $x \in S$ this formula gives the   BIE for problems solving ,
it is singular for the second BVP.

Using correlations (\ref{(3.4)})-(\ref{(3.6)}) the dynamic
analogue of Gauss  formula  also possible to write in integral
form.

L e m m a 6.1.  \textit{By N=3 }
\[
\int\limits_{S_t \left( x \right)} {\frac{{\partial r^{ - 1} }}
{{\partial n\left( y \right)}}dS\left( y \right)}  + \frac{1}
{c}\frac{\partial } {{\partial t}}\left\{ {\int\limits_{r = ct}
{\frac{{H_S^ -  \left( y \right)}} {r}dS\left( y \right) + }
\int\limits_{S_t \left( x \right)} {\frac{{\partial \ln r}}
{{\partial n\left( y \right)}}dS\left( y \right)} } \right\} =
4\pi H_S^ -  \left( x \right)H\left( t \right)
\]
\textit{When $t > t^* \left( x \right)$ thence  the known Gauss
formula }[ 2] follows:
\begin{equation}\label{(6.8)}
\int\limits_S {\frac{{\partial  }} {{\partial n( y )}}}
\left(\frac{1}{r}\right)\,\,dS\left( y \right) = 4\pi H_S^ -
\left( x \right).
\end{equation}

Proof of this formulas similarly is like to proof of the lemma
5.1. Notice  that for $N=2$  two-dimensional analogue of the Gauss
formula is not followed from lemmas 5.1.\vspace{3mm}

\textbf{7. Solution of the BVP for Dalamber equation  (N=1).}

T h e o r e m  7.1. \textit{ For $N=1$ BVP solution has the
following integral presentation}
\begin{equation}\label{(7.1)}
\begin{gathered}
2\hat u = c\, H\left( {ct - \left| {x - a_2 } \right|}
\right)\int\limits_{{{\left| {x - a_2 } \right|}/c}}^t {u,_x
\left( {a_2 ,\tau } \right)d\tau }  - c\,H\left( {ct - \left| {x -
a_1 } \right|} \right)\int\limits_{{{\left| {x - a_1 }
\right|}/c}}^t {u,_x \left( {a_1 ,\tau } \right)d\tau }  +
\\
+ \operatorname{sgn} \left({x - a_1 } \right)H\left( {ct - \left|
{x - a_1 } \right|} \right)u\left( {a_1 ,t - \frac{{\left| {x -
a_1 } \right|}} {c}} \right)-\\ - \operatorname{sgn} \left( {x -
a_2 } \right)H\left( {ct - \left| {x - a_2 } \right|}
\right)u\left( {a_2 ,t - \frac{{\left| {x - a_2 } \right|}} {c}}
\right) +
\\
+  c^{ - 1} \int\limits_{a_1 }^{a_2 } {\dot u_0 \left( y
\right)H(ct - \left| {x - y} \right|)dy}  + u_0 \left( {x + ct}
\right)H_S^ -  \left( {x + ct} \right) + u_0 \left( {x - ct}
\right)H_S^ -  \left( {x - ct} \right).
\end{gathered}
\end{equation}

P r o o f.  Denote $x_1 = x$. In this case ([2], p.206)
\begin{equation}\label{(7.2)}
\begin{gathered}
\hat U\left( {x,t} \right) =  - \frac{c} {2}H\left( {ct - \left| x
\right|} \right) ,\quad \hat U,_t  =  - \frac{c} {2}\delta \left(
{t - \left| x \right|/c} \right),
\\
\hat W\left( {x,t} \right) =  - \frac{c} {2}H\left( {ct - \left| x
\right|} \right)\left( {ct - \left| x \right|} \right), \quad \hat
W,_x \left( {x,t} \right) = \frac{c} {2}H\left( {ct - \left| x
\right|} \right)\operatorname{sgn} \,x ,
\end{gathered}
\end{equation}
\begin{equation}\label{(7.3)}
{\mathop{\rm sgn}} \,x = \left\{ \begin{array}{l}
 \;\,1,\quad x > 0; \\
 \,\,0,\quad x = 0; \\
  - 1,\;\,\,\,x < 0. \\
\end{array} \right.
\end{equation}
In this case  it is impossible to use the formula (\ref{(3.7)})
for building of the integral analogue of the Green formula ,
because in it some functions are not determined .

It is possible to get the similar formula, if put $u$ as zero
outside of given interval and considering an action of wave
operator on it in  the class of generalized functions. We enter
otherwise to use the formula (\ref{(3.5)}). We  increase the
 domain of definition $u\left( {x,t} \right)$ in the band in $R^2
\times R^ + $: $\{a_1 \leqslant \;x_1 \leqslant a_2$ , $ - \infty
< \;x_2 < \infty , t>0\}$. Then boundary$S$ will consist of two
direct lines $x_1 = a_1 ,\, x_1 = a_2 $, which external normals
have coordinates (-1,0) and (1,0) accordingly, $\frac{{\partial
u}}{{\partial n}} = n_1 \frac{{\partial u}}{{\partial x_1 }}$ for
$x \in S$, $H_S^ - \left( x \right) = H\left( {x_1  - a_1 }
\right)H\left( {a_2  - x_1 } \right)$, $n_1 \delta _S \left( x
\right) = \frac{{\partial H_S^ -  }}{{\partial x_1 }} =  - \delta
\left( {x_1  - a_1 } \right) + \delta \left( {x_1 - a_2 }
\right)$. The formula (\ref{(3.7)}) of the theorem 3.1 is
converted to type:
\begin{equation}\label{(7.4)}
\begin{gathered}
\hat u = \hat U_2  * \frac{{\partial u}} {{\partial x}}H\left( t
\right)\left( {\delta \left( {x - a_2 } \right) - \delta \left( {x
- a_1 } \right)} \right) + \frac{{\partial \hat W_2 }} {{\partial
x}} * \dot u\left( {x,t} \right)H\left( t \right)\left( {\delta
\left( {x - a_2 } \right) - \delta \left( {x - a_1 } \right)}
\right) +
\\+ c^{ - 2} \left( {\hat U_2 ,_t \mathop  * \limits_x u_0 \left( x
\right)H_S^ -  \left( x \right)} \right) + c^{ - 2} \hat U_2  *
\hat G.
\end{gathered}
\end{equation}
Here, all convolutions are taken with Green function for $N=2$ and
its antiderivative. On base of the method of the lowering on $x_2
$, rolling up on $x_2$, since $u$ does not depend on $x_2$ , we
get
\begin{equation}\label{(7.5)}
\begin{gathered}
\hat u = \hat U\left( {x - a_2 ,t} \right)\mathop  * \limits_t
\frac{{\partial u\left( {a_2 ,t} \right)}} {{\partial x}}H\left( t
\right) - \hat U\left( {x - a_1 ,t} \right)\mathop  * \limits_t
\frac{{\partial u\left( {a_1 ,t} \right)}} {{\partial x}}H\left( t
\right) +
\\
+\hat W,_x \left( {x - a_2 ,t} \right)\mathop  * \limits_t \dot
u\left( {a_2 ,t} \right)H\left( t \right) - \hat W,_x \left( {x -
a_1 ,t} \right)\mathop  * \limits_t \dot u\left( {a_1 ,t}
\right)H\left( t \right) + \hat W,_x \left( {x - a_2 ,t}
\right)u_0 \left( {a_2 } \right) - \\
-\hat W,_x \left( {x - a_1 ,t} \right)u_0 \left( {a_1 } \right) -
c^{ - 2} \left( {\hat U\mathop  * \limits_x \dot u_0 \left( x
\right)H_S^ -  \left( x \right)} \right) - c^{ - 2} \hat U,_t
\mathop  * \limits_x u_0 \left( x \right)H_S^ -  \left( x \right)
+ \hat U * \hat G.\\
\end{gathered}
\end{equation}
Substituting (\ref{(7.3)}) into (\ref{(7.5)}) and executing
integration,  we get formula (\ref{(7.1)}).

Easy to show that it is equitable also for $x = a_1 ,\,\,x = a_2 $
(with regard for (\ref{(7.3)})). For this it is enough to write
formula (\ref{(7.1)}) for interval $ \left( {a_1 ',a_2 } \right) =
\left( {a_1  + \varepsilon ,a_2 } \right)\, ( \left( {a_1 ,a_2 -
\varepsilon } \right)$). Supposing $x = a_1 $ and  $\varepsilon
\to + 0$ we have:
\[
\begin{gathered}
0 = \mathop {\lim }\limits_{\varepsilon  \to 0} c\left\{ {H\left(
{ct - d} \right)\int\limits_{\frac{{\left| {x - a_2 } \right|}}
{c}}^t {u,_x \left( {a_2 ,\tau } \right)d\tau }  - H\left( {ct -
\varepsilon } \right)\int\limits_{\frac{\varepsilon } {c}}^t {u,_x
\left( {a_1 ,\tau } \right)d\tau } } \right\} +\\+ H\left( {ct -
d} \right)u\left( {a_2 ,t - d/c} \right) + \operatorname{sgn}
\left( { - \varepsilon } \right)H\left( {ct} \right)u\left( {a_1
,t} \right) + c^{ - 1} \int\limits_{a_1  +
 \varepsilon }^{a_2 } {\dot u_0 \left( y \right)H(ct -
 \left| {a_1  - y} \right|)dy + }  \\
+ u_0 \left( {a_1  + ct} \right)H_S^ -  \left( {a_1  + ct} \right)
+ u_0 \left( {a_1  - ct} \right)H_S^ -  \left( {a_1  - ct} \right)
=
\\
= c\left\{ {H\left( {t - d/c} \right)\int\limits_{\frac{{\left| {x
- a_2 } \right|}} {c}}^t {u,_x \left( {a_2 ,\tau } \right)d\tau }
- \int\limits_0^t {u,_x \left( {a_1 ,\tau } \right)d\tau } }
\right\} + H\left( {t - d/c} \right)u\left( {a_2 ,t - \frac{d}
{c}} \right) +  \\
- H\left( t \right)u\left( {a_1 ,t} \right) + c^{ - 1}
\int\limits_{a_1 }^{a_2 } {\dot u_0 \left( y \right)H(t - \left|
{a_1  - y} \right|/c)dy}  + u_0
\left( {a_1  + ct} \right)H_S^ -  \left( {a_1  + ct} \right) +  \\
+ u_0 \left( {a_1  - ct} \right)H_S^ -  \left( {a_1  - ct} \right) \\
\end{gathered}
\]
($d = \left| {a_1 - a_2 } \right|$). Transferring the  summand $ -
H\left( t \right)u\left( {a_1 ,t} \right)$ into the left part,
with regard for values of the characteristic function on border,
we get the formula of the theorem for left endpoint. By the
similar way the formula is proved for $x = a_2 $. As a result on
the end of the interval$(a_1,a_2)$ we have the following equations
for determination of unknown boundary functions:
\[
u(a_1 ,t) = c {H\left( {ct - d} \right)\int\limits_{d/ c}^t
{u_{,x} \left( {a_2 ,\tau } \right)d\tau }  - cH(t)\int\limits_0^t
{u_{,x} \left( {a_1 ,\tau } \right)d\tau } }  +
 H\left( {ct - d} \right)u\left( {a_2 ,t - \frac{d}
{c}} \right) +\]
\[ +c^{ - 1} \int\limits_{a_1 }^{a_2 } {\dot u_0
\left( y \right)H(ct - \left| {a_1  - y} \right|)dy}  + u_0 \left(
{a_1  + ct} \right)H\left( {d - ct} \right)H\left( t \right)\quad
\textrm{для} \,\,x = a_1;
\]
\[
u(a_2 t) = c {H\left( t \right)\int\limits_0^t {u,_x \left( {a_2
,\tau } \right)d\tau }  - cH\left( {ct - d}
\right)\int\limits_{{d/c}}^t {u,_x \left( {a_1 ,\tau }
\right)d\tau } }  + H\left( {ct - d} \right)u\left( {a_1 ,t -
\frac{d}{c}} \right) +
\]
\[
+c^{ - 1} \int\limits_{a_1 }^{a_2 } {\dot u_0 \left( y \right)H(ct
- \left| {a_2 - y} \right|)dy} + u_0 \left( {a_2 - ct}
\right)H\left( {d - ct} \right)H\left( {ct} \right)\quad
\textrm{для} \,\, x = a_2.
\]

Under given $u,_x \left( {a_k ,t} \right),\,k = 1,2$, we get two
functional equations with lagging argument for determination $u$
on boundary of the area, which can be solved incremental on time
from $t=0$. Under the known $u\left( {a_1 ,t}
\right),\,\,\,u\left( {a_2 ,t} \right)$ we have  the system of two
integral equations.  \vspace{3mm}

\textbf{Conclusion.}   Using this method it is possible to build
the similar formulas and boundary integral equations for solutions
of BVP in space of greater dimension( $N>3$). It is particularly
efficient  in BVP for systems of the equations of mathematical
physics, when  the Green matrix of system can me constructed.
Herewith the type of the equations is unessential , it can be also
elliptical, as, for instance in problem of stationary diffraction
of  electromagnetic waves [5], or in problem of the  elasticity
theory [6], parabolic or mixed type in problem of
thermoelastodynamics [7]. But particularly efficient this method
for solving of the hyperbolic equations, where using classical
methods more difficult, but sometimes and simply impossible (refer
to, for instance, [6,8,9]).

The study of solubility of built BIEs presents the independent
problems of the functional analysis since these equations do not
pertain to well studied classical ones. However we notice that
using the computing methods on base of the methods of  boundary
element with transition to discrete analogue BIE, allows
effectively to build the solutions of like problems [10].
\vspace{4mm}

\centerline{\textbf{References}}

\small { [1] Petrovskiy I.S. The lectures about equations with
partial derivatives. M., 1961

[2]  Vladimirov V.S. The equations of mathematical physics. M.,
1978, 512 p.

[3] Alekseeva L.A. Integral equations of  boundary value problems
for wave equation in $R_2\times t$ // Differential equations,
1992,V.28,\textbf{8}.

[4]  Alekseeva L.A. The dynamic analogues of Green and Gauss
formulas for solutions of the wave equation in $R_N\times t$
//Differential equations.1995. V.31,  \textbf{11}.

[5] Alexeyeva L.A., Sautbekov S.S.  The method of generalized
function at solving of the stationary boundary value problems for
Maxwell equations // Journal of computational mathematics and
mathematical physics. 2000.Т.40, \textbf{ 4}.

[6]  Alexeyeva L.A. Boundary Element Method of Boundary Value
Problems of  elasto\-dy\-na\-mics by stati\-ona\-ry run\-ning
loads//Int. J. Engineering Analysis with Boundary Element. 1998,
\textbf{11}.

[7]  Alekseeva L.A., Kupesova B.N. The method of generalized
function in boundary value problems of thermoelastodynamics
// Applied mathematics and mechanics. 2001. V.65,\textbf{2}.

[8]  Alexeyeva L.A. The Generalized solutions of nonstationary
boundary value problems for Maxwell equations  // Journal of
computational mathematics and mathematical physics. 2002.  V.42,
\textbf{1}.

[9] Alexeyeva L.A., Zakiryanova G.K. Generalized solutions of
boundary value problems of dynamics of anisotropic  elastic
media// Journal of the Mechanical Behavior of Materials. 2004.
 \textbf{5}.

[10] Alexeyeva L.A., Dildabaev Sh.A., Zhanbyrbaev A.B.,
Zakiryanova G.K. Boundary Integral Equation Method in Two and
three dimensional problems of elastodynamics// Int. J.
Computational Mechanics.1996.V.18,\textbf{2}. }

\vspace{10mm}

Lyudmila A. Alexeyeva\\
Professor, head of waves dynamics laboratory,\\
Institute of Mathematics of Education and Sciences Ministry,\\
Pushkin str. 125, Almaty, 050010, Kazakhstan,\\

tel. +7 3272 911624,

E-mail: alexeeva@math.kz

\vspace{5mm}

\newpage

\textbf{Alexeyeva L.A.  Nonstationary  boundary value problems for
wave equation and their generalized solutions}
\vspace{3mm}

\textbf{Abstract}. The multivariate analogue of Dalamber's equation  in the space of
generalized functions is considered. The method of
generalized functions for the building of solutions of
nonstationary boundary value problems for wave equations in spaces
of different dimensions is elaborated.  Dynamic analogues of
Green and Gauss formulas for solutions of  wave equation in the
space of generalized functions are built.  Their regular integral
repre\-sen\-ta\-tions and singular  boundary integral equations for
solving the nonstationary problems are constructed for the spaces of the dimensions 1,2,3.
The method of obtaining  of conditions on fronts of shock waves is stated.

\end{document}